\documentclass[12pt,reqno,a4paper]{amsart}
\usepackage{amssymb,amsfonts,latexsym,mathrsfs}

\title[Coxeter groups, imaginary cones and dominance]{Coxeter groups, imaginary cones and dominance}

\author{Fu, Xiang}

\dedicatory{\upshape
School of Mathematics and Statistics\\
University of Sydney, NSW 2006, Australia\\[.5em]
\texttt{xifu9119@mail.usyd.edu.au}\\
\texttt{X.Fu@maths.usyd.edu.au}\\[1em]
Preliminary version,
\today
}

\newtheorem{theorem}{Theorem}[section]
\newtheorem{lemma}[theorem]{Lemma}
\newtheorem{proposition}[theorem]{Proposition}
\newtheorem{corollary}[theorem]{Corollary}

\theoremstyle{definition}
\newtheorem{definition}[theorem]{Definition}
\newtheorem{notations}[theorem]{Notations}

\theoremstyle{remark}
\newtheorem{remark}[theorem]{Remark}

\numberwithin{equation}{section}

\newcommand{\Z}{\mathbb{Z}}
\newcommand{\N}{\mathbb{N}}
\newcommand{\R}{\mathbb{R}}

\DeclareMathOperator{\dom}{\, dom}
\DeclareMathOperator{\PLC}{PLC}
\DeclareMathOperator{\coeff}{coeff}
\DeclareMathOperator{\pos}{Pos}
\DeclareMathOperator{\n}{Neg}
\DeclareMathOperator{\supp}{supp}
\DeclareMathOperator{\GL}{GL}
\DeclareMathOperator{\spa}{span}
\DeclareMathOperator{\Hom}{Hom}
\DeclareMathOperator{\dep}{dp}


\begin{document}

\thanks{The work presented in this paper was completed when the author was supported by 
an Australian Research Council Discovery Project: \emph{Invariant theory, cellularity and geometry}, 
No. DP0772870.}

\begin{abstract}
Brink and Howlett have introduced a partial ordering, called \emph{dominance},
on the root systems of Coxeter groups in their proof that all finitely 
generated Coxeter groups are automatic (Math. Ann. \textbf{296} (1993), 179--190).
Recently a function called
\emph{$\infty$-height} is defined on the reflections of  Coxeter
groups in an investigation of various regularity 
properties of Coxeter groups
(Edgar, \emph{Dominance and regularity in Coxeter groups}, PhD thesis, 2009).  
In this paper, we show that these two concepts are closely related to each other.
We also give applications of
dominance to the study of \emph{imaginary cones} of Coxeter groups.
\end{abstract}

\maketitle

\section{Introduction}
In this paper we attempt to extend the understanding of a partial ordering 
(called \emph{dominance}) defined on the root system of an arbitrary Coxeter
group. The dominance ordering was introduced by Brink and Howlett in their 
paper \cite{BH93} (where it was used to prove the automaticity of all finitely generated
Coxeter groups). Dominance ordering has been further studied in the 1990's by Brink (\cite{BB98})
and Krammer (\cite{DK94}, and later reproduced in \cite{DK09}), and it has only been 
recently examined again (Dyer \cite{MD12},
in connection with the representation theory of Coxeter groups; 
the PhD thesis of Edgar \cite{TE09}; and a recent paper by the 
author \cite{FU1}). The present paper is a short addition to both \cite{TE09} and \cite{FU1}, 
and it could serve as a building block in the general knowledge on dominance ordering and on 
the combinatorics and geometry of Coxeter groups in general.

More specifically, this paper has the following two objectives: (1) investigating the connection 
between the dominance ordering on the root system of an arbitrary Coxeter groups $W$ and a 
specific function (called \emph{$\infty$-height}) defined on the set of reflections of $W$; (2) 
exploring the applications of the dominance ordering to the \emph{imaginary cone} of $W$ 
(as defined by Kac).

The paper is organized into three sections. In the first section, background material 
is introduced: root basis, Coxeter datum, and root systems are defined in the context of
the paper, and some basic properties of Coxeter groups are recalled for later use in the paper
(most of them can be found in Howlett's lectures \cite{RB96}). Here we follow the definition used in \cite{DK94}, 
which gives a slight variant of the classical notion of root system, 
particularly adapted when working with arbitrary (not necessarily crystallographic) 
Coxeter groups. Furthermore, this framework allows easy passing 
to reflection subgroups. Indeed, we recall the fundamental property (\cite[Theorem 1.8]{MD87}) 
that the reflection subgroups of a Coxeter group are themselves Coxeter groups, and this
particular framework allows us to apply all the definitions and properties to the reflection 
subgroups and not only to the over-group.

In the second section, the first main theorem (giving the connection between $\infty$-height 
and dominance order) is stated and proved. All results are related to an arbitrary 
Coxeter datum, implying the data of a root system $\Phi$, its associated Coxeter group $W$,
and the set $T$ of all reflections of $W$ (consisting of all the $W$-conjugates of the Coxeter
generators). The main objects of study are:
\begin{itemize}
\item{} the dominance order on $\Phi$ (Definition \ref{def:dom}): given
$x,y\in \Phi$, we say $x$ \emph{dominates} $y$ if whenever $w\in W$ such that
$wx\in\Phi^-$ then $wy\in \Phi^-$ too (where $\Phi^-$ denotes the set of negative
roots);
\item{} the function $\infty$-height on $T$. It is a variant of the usual (standard) height 
function of a reflection $t\in T$, namely, the minimal length of an element of $W$ 
that maps $\alpha_t$ (the unique positive root associated to $t$) to an element of 
the root basis. Adhering to the general framework of this paper, our definition of 
the height function applies to all reflection subgroups of $W$.
It is easy to check (Lemma ~\ref{lem:height}) that the height of $t$ is equal to 
the sum of the heights of $t$ relative to each maximal (with respect to inclusion) 
dihedral reflection subgroup containing $t$. The $\infty$-height of $t$ is then defined
as a sub-sum of this sum, taking into account only those subgroups 
which are infinite (Definition~\ref{def:h}). 
\end{itemize}
We then show that these two concepts are closely related in the following way. The 
canonical bijection $t\leftrightarrow \alpha_t$, between $T$ and $\Phi^+$ (the set of 
positive roots), restricts to a bijection between (for any $n\in \N$):
\begin{itemize}
\item{} the set $T_n$ of all reflections whose $\infty$-height is $n$; and
\item{} the set $D_n$ of all positive roots which strictly dominate exactly $n$ other 
positive roots.
\end{itemize}
The proof of this fact (Theorem \ref{th:bij}) relies on a study of dihedral
reflection subgroups. We have previously studied the partition $(D_n)_{n\in \N}$
of $\Phi^+$ in \cite{FU1}; in particular, we showed there that each $D_n$ is finite
and we gave an upper bound for its cardinality. Together with Theorem~\ref{th:bij}, 
this allows us to deduce here some information on the combinatorics of the $T_n$'s
(Corollary \ref{co:T_n}).

The final section explores the relation between the dominance order and
the imaginary cone of a Coxeter group. The concept of \emph{imaginary cone}
was introduced by Kac in \cite{VK} to study the imaginary roots of Kac-Moody
Lie algebras, and was later generalized to Coxeter groups by H\'ee \cite{HE1,HE2} 
and Dyer \cite{MD12}. It is defined as the subset of the dual of the Tits cone
(denoted as $U^*$ here)  consisting of elements $v\in U^*$ such that 
$(v, \alpha)>0$ for only finitely-many $\alpha\in \Phi^+$ (where $(\,,\,)$ denotes
the bilinear form associated to the Coxeter datum). The main results (Theorem~\ref{th:imc} and 
Corollary~\ref{co:imc})
of this section state the following property: whenever $x, y\in \Phi$, then $x$ dominates $y$ 
if and only if $x-y$ lies in the imaginary cone. One direction of this property was first suggested to us
by Howlett (private communications), and it is a special case of a result 
obtained independently (but earlier) by Dyer. We are deeply indebted to both of them for helpful 
discussions inspiring us to study the 
imaginary cone. We would also like to thank the referee of this paper
for many valuable suggestions, especially those resulting in Corollary~\ref{co:imc}. 
  To close this section, we include an alternative definition for the imaginary cone 
in the case where $W$ is finitely generated.

\section{Background Material}
\begin{definition}\textup{(Krammer \cite{DK94})}
 \label{def:datum}
Suppose that $V$ is a vector space over $\R$ and let $(\,,\,)$ be a bilinear
form on $V$ and let $\Delta$ be a subset of $V$. Then $\Delta$ is called a \emph{root
basis} if the following conditions are satisfied:
\begin{itemize}
 \item [(C1)] $(a, a)=1$ for all $a\in \Delta$, and for distinct elements $a, b\in\Delta$
  either $(a, b)=-\cos(\pi/m_{ab})$ for some integer
$m_{ab}=m_{ba}\geq 2$, or else $(a, b) \leq -1$ (in which case we define
$m_{ab}=m_{ba}=\infty$);
 \item [(C2)] $0\notin \PLC(\Pi)$, where $\PLC(A)$, \emph{the positive linear cone of}
 a set $A$,  denotes the set
 $$\{\,\sum\limits_{a\in A} \lambda_a a\mid \text{$\lambda_a \geq 0$ for all
$a\in A$ and $\lambda_{a'}>0$ for some $a'\in A$}\,\}.$$ 
\end{itemize}
\end{definition}

If $\Delta$ is a root basis, then we call the triple $\mathscr{C}=(\,V, \, \Delta,
\,(\,,\,)\,)$ a \emph{Coxeter datum}. Throughout this paper we fix a particular
Coxeter datum $\mathscr{C}$.  We stress that our definition of a root basis is not the most
classical one of \cite{BN68} or even \cite{HM}: the root system (see Definition~\ref{def:rootsys}) 
arising from our definition of a root basis is not necessarily crystallographic (indeed, the
bilinear form can take values less than $-1$), and the root basis is not assumed to be linearly
independent (this allows us to transmit easily the definitions and properties of a Coxeter group 
to its reflection subgroups, indeed the requirements in our definition of a root basis of a Coxeter
group are identical to those in the characterization of the equivalent of a root basis in any
reflection subgroup).
Observe that (C1) implies that for each $a\in \Delta$, $a\notin
\PLC(\Delta\setminus\{a\})$, and furthermore, (C1) together with (C2) yield that $\{a,
b, c\}$ is linearly independent for all distinct $a, b, c\in \Delta$. 
Note also that (C2) is equivalent to the requirement that 
$0$ does not lie in the convex hull of $\Delta$.

For each $a\in \Delta$,
define $\rho_a \in \GL(V)$ by the rule: $\rho_a x=x-2(x, a)a$, for all $x\in V$.
Observe that $\rho_a$ is a reflection, and $\rho_a a=-a$. The following
proposition summarizes a few useful results:
 
\begin{proposition}\textup{\cite[Lecture 1]{RB96}}
 \label{pp:anu1}
\rm{(i)}\quad Suppose that $a, b\in \Delta$ are distinct such that $m_{ab}\neq
\infty$. Set $\theta =\frac{\pi}{m_{ab}}$. Then 
$$(\rho_a \rho_b)^i a=\frac{\sin(2i+1)\theta}{\sin \theta}a+\frac{\sin
2i\theta}{\sin\theta}b, $$
for each integer $i$, and in particular, $\rho_a \rho_b$ has order $m_{ab}$ in
$\GL(V)$.

\noindent\rm{(ii)}\quad Suppose that $a, b\in \Delta$ are distinct such that
$m_{ab}=\infty$. Set $\theta =\cosh^{-1}(-(a, b))$. Then  
\begin{equation*}
(\rho_a \rho_b)^i a=
\begin{cases}
\frac{\sinh(2i+1)\theta}{\sinh \theta}a+\frac{\sinh
2i\theta}{\sinh\theta}b, \text{ if $\theta \neq 0$}\\
(2i+1) a+2i b, \text{ }\text{ }\text{ }\text{ }\text{ }\text{ }\text{ }\text{ }\text{if $\theta =0$},
\end{cases}
\end{equation*}
for each integer $i$, and in particular, $\rho_a \rho_b$ has infinite order in
$\GL(V)$.
\qed
\end{proposition}

Let $G_{\mathscr{C}}$ be the subgroup of $\GL(V)$ generated by 
$\{\,\rho_a\mid a\in \Delta\,\}$.
Suppose that $(W, S)$ is a Coxeter system in the sense of \cite{HH81} or \cite{HM} with
$S=\{\,r_a\mid a\in \Delta\,\}$ being a set of involutions generating $W$ 
subject only to the condition that
the order of $r_a r_b$ is $m_{ab}$ for all $a, b\in \Delta$ with $m_{ab}\neq \infty$. Then
Proposition~\ref{pp:anu1} yields that there exists a group homomorphism
$\phi_{\mathscr{C}}\colon W\to G_{\mathscr{C}}$ satisfying
$\phi_{\mathscr{C}}(r_a)=\rho_a$ for all $a\in \Delta$. This homomorphism together
with the $G_{\mathscr{C}}$-action on $V$ give rise to a $W$-action on $V$: for
each $w\in W$ and $x\in V$, define $wx\in V$ by $wx=\phi_{\mathscr{C}}(w)x$. It
can be easily checked that this $W$-action preserves $(\,,\,)$.
Denote the length function of $W$ with respect to $S$ by $\ell$, and call an 
expression $w=r_1 r_2\cdots r_n$ (where $w\in W$ and $r_i\in S$) \emph{reduced} if
$\ell(w)=n$. The following is a useful result:

\begin{proposition}\textup{\cite[Lecture 1, Theorem, Page 4]{RB96}}
 \label{pp:anu2}
Let $G_{\mathscr{C}}, W, S$ and $\ell$ be as  above, and let $w\in W$ and $a\in \Delta$. If
$\ell(wr_a)\geq \ell(w)$ then $wa\in \PLC(\Delta)$.
\qed
\end{proposition}
An immediate consequence of the above proposition is the following important fact:
\begin{corollary}\textup{\cite[Lecture 1, Corollary, Page 5]{RB96}}
 \label{co:anu2}
Let $G_{\mathscr{C}}, W, S$ and $\phi_{\mathscr{C}}$ be as above. Then
$\phi_{\mathscr{C}}\colon W\to G_{\mathscr{C}}$ is an isomorphism.
\qed
\end{corollary}

In particular, the above corollary yields that $(G_{\mathscr{C}}, \{\,\rho_a\mid
a\in \Delta\, \})$ is a Coxeter system isomorphic to $(W, S)$. We call $(W, S)$ the
\emph{abstract Coxeter system} associated to the Coxeter datum $\mathscr{C}$, and
we call $W$ a \emph{Coxeter group} of rank $\#S$ (where $\#$ denotes
cardinality).

\begin{definition}
\label{def:rootsys}
The \emph{root system} of $W$ in $V$ is the set 
$$\Phi=\{\,wa \mid \text{$w\in W$ and $a\in \Delta$}\,\}.$$
The set $\Phi^+=\Phi\cap \PLC(\Delta)$ is called the set of \emph{positive roots},
and the set $\Phi^-=-\Phi^+$ is called  the set of \emph{negative roots}.
\end{definition}
From Proposition \ref{pp:anu2} we may readily deduce that:
\begin{proposition}\textup{(\cite[Lecture 3]{RB96})}
\label{pp:anu3}
\rm{(i)}\quad Let $w\in W$ and $a\in \Delta$. Then 
 \begin{equation*}
\ell(wr_a) =
\begin{cases}
\ell(w)-1  \text{,   if } wa\in \Phi^-,\\
\ell(w)+1  \text{,   if } wa\in \Phi^+.
\end{cases}
\end{equation*}

\noindent\rm{(ii)}\quad $\Phi=\Phi^+\biguplus\Phi^-$, where $\biguplus$ denotes
disjoint union.

\noindent\rm{(iii)}\quad $W$ is finite if and only if $\Phi$ is finite. 
\qed
\end{proposition}

Define $T=\bigcup_{w\in W} w S w^{-1}$. We call $T$ the set of \emph{reflections}
in $W$. For each $x\in \Phi$, let $\rho_x\in \GL(V)$ be defined by the rule: $\rho_x
(v)=v-2(v,x)x$, for all $v\in V$. Since $x\in \Phi$, it follows that $x=wa$ for
some $w\in W$ and $a\in \Delta$. Direct calculations yield that $\rho_x
=(\phi_{\mathscr{C}} (w)) \rho_a (\phi_{\mathscr{C}} (w))^{-1}\in G_{\mathscr{C}}$.
Now let $r_x\in W$ be such that $\phi_{\mathscr{C}}(r_x)=\rho_x$. Then $r_x = w r_a
w^{-1}\in T$ and we call it the reflection corresponding to $x$. It is readily
checked that $r_x =r_{-x}$ for all $x\in \Phi$ and $T=\{r_x \mid x\in \Phi\}$.
For each $t\in T$ we let $\alpha_t$ be the unique positive root with the
property that $r_{\alpha_t}=t$. It is also easily checked that there is a
bijection $\psi\colon T\to \Phi^+$ given by $\psi(t)= \alpha_t$, and we call
$\psi$ the \emph{canonical bijection}.

For each $x\in \Phi^+$, as in \cite{BH93}, we define the \emph{depth} of $x$
relative to $S$ to
be $\min\{\, \ell(w) \mid\text{$w\in W$ and $wx\in \Phi^-$}\, \}$, and we denote
it by $\dep(x)$. The following lemma gives some basic properties of depth:
\begin{lemma}\textup{(\cite{BH93, BB94, SC91}).}
\label{lem:pre}
\begin{itemize}
 \item[(i)] Let $\alpha\in \Phi^+$. Then $\dep(\alpha)=\frac{1}{2}(\ell(r_{\alpha})+1)$.
 \item[(ii)] Let $r\in S$ and $\alpha \in \Phi^+\setminus\{\alpha_r\}$. Then
\begin{equation*}
\dep(r\alpha) =
\begin{cases}
\dep(\alpha)-1  \text{   if } (\alpha, \alpha_r)>0,\\
\dep(\alpha) \ \ \ \ \ \text{         if } (\alpha, \alpha_r)=0,\\
\dep(\alpha)+1 \text{   if } (\alpha, \alpha_r)<0.
\end{cases}
\end{equation*}
\end{itemize}
\end{lemma}
\begin{proof}
\rm{(ii)}:  \cite[Corollary 2.7]{BB94}.

\noindent\rm{(ii)}: \cite[Lemma 1.7]{BH93}.
\end{proof}

\begin{remark}
Part (i) of the above Lemma is equivalent to the property that any reflection in a Coxeter group 
has a palindromic expression which is reduced, and this was indeed noted in \cite[Proposition 4.3]{SC91}. 
\end{remark}

Define functions $N\colon W\to \mathcal{P}(\Phi^+)$ and
$\overline{N}\colon W\to \mathcal{P}(T)$ (where $\mathcal{P}$ denotes power set)
by setting
$N(w)=\{\,x\in \Phi^+\mid wx\in \Phi^-\,\}$ and
$\overline{N}(w)=\{\,t\in T \mid \ell(wt)<\ell(w)\,\}$ for all $w\in W$. We call
$\overline{N}$ the \emph{reflection cocycle} of $W$ (sometimes $\overline{N}(w)$ is also called
the \emph{right descent set} of $w$). Standard arguments as those
in \cite[$\S$~5.6]{HM} yield that for each $w\in W$, 
\begin{align}
 \label{eq:ell1}
\ell(w)&=\#N(w),\\
\noalign{\hbox{and}}
\overline{N}(w)&=\{\,r_x\mid x\in N(w)\,\}.
\end{align}
In particular, $N(r_a)=\{a\}$ for
$a\in \Delta$. Moreover, $\ell(wv^{-1})+\ell(v)=\ell(w)$, for some $w, v\in W$ if
and only if $N(v)\subseteq N(w)$.
  
A subgroup $W'$ of $W$ is a \emph{reflection subgroup} of $W$ if  $W'=\langle
W'\cap T\rangle$ ($W'$ is generated by the reflections contained in it). For
any reflection subgroup $W'$ of $W$, let 
\begin{align*}
S(W')&=\{\,t\in T\mid \overline{N}(t)\cap W'=\{t\}\,\}\\
\noalign{\hbox{and}}
\Delta(W')&=\{\,x\in \Phi^+\mid r_x\in S(W')\,\}.
\end{align*}
It was shown by Dyer (\cite{MD90}) and Deodhar (\cite{VD82}) that $(W', S(W'))$
forms a Coxeter system:

\begin{theorem}\textup{(Dyer)}
\label{th:croots}
\rm{(i)}\quad Suppose that $W'$ is an arbitrary reflection subgroup of $W$. Then 
 $(W', S(W'))$ forms a Coxeter system. Moreover,
$W'\cap T=\bigcup_{w\in W'}w S(W') w^{-1}$.

\noindent\rm{(ii)} Suppose that $W'$ is a reflection subgroup of $W$, and suppose
that $a, b\in \Delta(W')$ are distinct. Then 
$$
(a, b)\in \{\,-\cos(\pi/n)\mid \text{$n\in \N$ and $n\geq 2$}\,\}\cup (-\infty,
-1].
$$
And conversely if $\Delta'$ is a subset of $\Phi^+$ satisfying the condition
that 
$$
(a, b)\in \{\,-\cos(\pi/n)\mid \text{$n\in \N$ and $n\geq 2$}\,\}\cup (-\infty,
-1]
$$
for all $a, b\in \Delta'$ with $a\neq b$, then $\Delta' =\Delta(W')$ for some
reflection subgroup $W'$ of $W$. In fact, $W'=\langle \{\, r_a\mid a\in
\Delta'\,\}\rangle$. 
\end{theorem}
\begin{proof}
\rm{(i)}\quad \cite[Theorem 3.3]{MD90}.

\noindent\rm{(ii)} \cite[Theorem 4.4]{MD90}.
\end{proof}

Let $(\,,\,)'$ be the restriction of $(\,,\,)$ on the subspace
$\spa(\Delta(W'))$.  Then $\mathscr{C'}=(\,\spa(\Delta(W')),\,
\Delta(W'),\, (\,,\,)'\,)$ is a Coxeter datum with $(W', S(W'))$ being the
associated abstract Coxeter system. Thus the notion of a root system applies to 
$\mathscr{C'}$. We let $\Phi(W')$, $\Phi^+(W')$ and $\Phi^-(W')$ be,
respectively, the set of roots, positive roots and negative roots for the datum
$\mathscr{C'}$. Then $\Phi(W')=W'\Delta(W')$ and Theorem~\ref{th:croots}~(i)
yields that $\Phi(W')=\{x\in \Phi\mid r_x \in W'\}$. Furthermore, we have
$\Phi^+(W')=\Phi(W')\cap \PLC(\Delta(W'))$ and $\Phi^-(W')=-\Phi^+(W')$. We
call $S(W')$ the set of \emph{canonical generators} of $W'$, and we call
$\Delta(W')$ the set of \emph{canonical roots} of $\Phi(W')$. In this paper a
reflection subgroup $W'$ is called a \emph{dihedral reflection subgroup} if
$\#S(W')=2$.

A subset $\Phi'$ of $\Phi$ is called a \emph{root subsystem} if $r_y x\in \Phi'$
whenever $x, y$ are both in $\Phi'$. It is easily seen that there is a bijective
correspondence between the set of reflection subgroups $W'$ of $W$ and the set of root subsystems
$\Phi'$ of $\Phi$: $W'$ uniquely determines the root subsystem $\Phi(W')$, and $\Phi'$ 
uniquely determines the reflection subgroup  $\langle\{\,
r_x\mid x\in \Phi'\,\}\rangle$. 

The notion of a length function also applies to the Coxeter system $(W',
S(W'))$, and we let $\ell_{(W',\ S(W'))}\colon W'\to \N$ be the length function
for $(W', S(W'))$. If $w\in W'$ and $a\in \Delta(W')$ then applying 
Proposition~\ref{pp:anu3} to the Coxeter datum $\mathscr{C}'=(\,\spa(\Delta(W'))$ yields
\begin{equation}
\label{eq:length}
\ell_{(W', \,S(W'))}(wr_a) =
\begin{cases}
\ell_{(W', \,S(W'))}(w)-1  \text{,   if } wa\in \Phi^-(W'),\\
\ell_{(W', \,S(W'))}(w)+1  \text{,   if } wa\in \Phi^+(W').
\end{cases}
\end{equation} 
Similarly the notion of a reflection cocycle also applies 
to the Coxeter system $(W', S(W'))$. Let 
$\overline{N}_{(W',\, S(W'))}\colon W\to \mathcal{P}(W'\cap T)$ 
denote the reflection cocycle for $(W', S(W'))$. Then for each $w\in W'$, 
$$
\overline{N}_{(W', \, S(W'))}(w)=\{\, t\in W'\cap T\mid \ell_{(W',\, S(W'))}(wt)<\ell_{(W',\, S(W'))}(w)\,\}.
$$
And we define $N_{(W',\,S(W'))}(w)=\{\,x\in \Phi^+(W')\mid wx\in \Phi^-(W') \,\}$, for each $w\in W'$.
It is shown in \cite{MD87} that 
$\overline{N}_{(W',\ S(W'))}(w)=\overline{N}(w)\cap W'$
for arbitrary reflection subgroup $W'$ of $W$. Furthermore, it is readily seen that the canonical
bijection $\psi$ restricts to a bijection $\psi'\colon T\cap W'\to \Phi^+(W')$
given by $\psi'(t)=\alpha_t$. For $w\in W'$, applying (\ref{eq:ell1}) to the
Coxeter datum $\mathscr{C}'=(\,\spa(\Delta(W'),\, \Delta(W'),\, (\,,\,)'\, )$ yields
that 
\begin{equation}
\label{eq:ell}
\ell_{(W',\, S(W'))}(w)=\# N_{(W', \,S(W'))}(w).
\end{equation}
Furthermore, $\ell_{(W',\,S(W'))}(wv^{-1})+\ell_{(W',\,S(W'))}(v)=\ell_{(W',\,S(W'))}(w)$, 
for some $w, v\in W'$, 
precisely when $N_{(W', \,S(W'))}(v)\subseteq N_{(W', \,S(W'))}(w)$.

For a Coxeter datum $\mathscr{C}=(\, V,\, \Delta, \, (\,,\,)\,)$, since $\Delta$ may
be linearly dependent, the expression of a root in $\Phi$ as a linear
combination of elements of $\Delta$  may not be unique. Thus the concept of the
coefficient of an element of $\Delta$ in any given root in $\Phi$  is potentially
ambiguous. We close this section by specifying a canonical way of expressing a
root in $\Phi$  as a linear combination of elements from $\Delta$. This canonical
expression follows from a standard construction similar to that considered 
in \cite[Proposition 2.9]{HT97}.

Given a Coxeter datum $\mathscr{C}=(\, V,\, \Delta, \, (\,,\,)\,)$, let $E$ be a
vector space over $\R$ with basis $\Delta_E=\{\,e_a\mid a\in \Delta\,\}$ in bijective
correspondence with $\Delta$, and let $(\,,\,)_E$ be the unique bilinear form on $E$
satisfying 
$$
(e_a, e_b)_E =(a, b) \text{ for all } a, b\in \Delta.
$$
Then $\mathscr{C}_E=(\, E,\, \Delta_E, \, (\,,\,)_E\,)$ is a Coxeter datum.
Moreover, $\mathscr{C}_E$ and $\mathscr{C}$ are
associated to the same abstract Coxeter system $(W, S)$; indeed Corollary~\ref{co:anu2}
yields that the abstract Coxeter group $W$ is isomorphic to both
$G_{\mathscr{C}}=\langle \{\, \rho_a\mid a\in \Delta\,\}\rangle$ and
$G_{\mathscr{C}_E}=\langle\{\, \rho_{e_a}\mid a\in \Delta\,\}\rangle$. Furthermore,
$W$ acts faithfully on $E$ via $r_a y =\rho_{e_a} y$ for all $a\in \Delta$ and
$y\in E$.

Let $f\colon E\to V$ be the unique linear map satisfying $f(e_a)= a$, for all
$a\in \Delta$. It is readily checked that $(f(x), f(y))=(x, y)_E$, for all $x, y\in
E$. 
Now  for all $a\in \Delta$ and $y\in E$, 
\begin{align*}
 r_a(f(y))=\rho_a(f(y))=f(y)-2(f(y), a)a &=f(y)-2(f(y), f(e_a))f(e_a)\\
                       &=f(y-2(y, e_a)_E e_a)\\
                       &=f(r_a y).
\end{align*}
Then it follows that $w(f(y))=f(wy)$, for all $w\in W$ and all $y\in E$, since
$W$ is generated by $\{\,r_a\mid a\in \Delta\,\}$.
Let $\Phi_E$ denote the root system associated to the datum $\mathscr{C}_E$.
Standard arguments yield that:

\begin{proposition}\textup{\cite[Proposition 2.1]{FU1}}
\label{pp:eqv}
The restriction of $f$ defines a $W$-equivariant bijection
$\Phi_E\leftrightarrow \Phi$.
\qed
\end{proposition}

Since $\Delta_E$ is linearly independent, it follows that each root $y\in \Phi_E$
can be written uniquely as $y=\sum_{e_a\in \Delta_E} \lambda_a e_a$; we say that
$\lambda_a$ is the \emph{coefficient} of $e_a$ in $y$, and it is denoted by
$\coeff_{e_a}(y)$. We use this fact together with the $W$-equivariant bijection
$f\colon \Phi_E\leftrightarrow \Phi$ to give a canonical expression of a root in
$\Phi$ in terms of $\Delta$:

\begin{definition}
 Suppose that $x\in \Phi$. For each $a\in \Delta$, define the \emph{canonical
coefficient} of $a$ in $x$, written $\coeff_a(x)$, by requiring that
$\coeff_a(x)=\coeff_{e_a}(f^{-1}(x))$. The \emph{support}, written $\supp(x)$, is
the set of $a\in \Delta$ with $\coeff_a(x)\neq 0$.
\end{definition}

\section{Dominance, Maximal Dihedral Reflection Subgroups and Infinity Height}
\label{sec:bij}

Throughout this section, let $W$ be the abstract Coxeter group associated to the Coxeter datum 
$\mathscr{C}=(\,V,\, \Delta, \,(\,,\,)\,)$, and let $\Phi$ and $T$ be the
corresponding root system and the set of reflections respectively. Recently in \cite{TE09},
a uniquely determined non-negative integer, called \emph{$\infty$-height}, 
is assigned to each reflection in $W$. Naturally, the set $T$ is then the 
disjoint union of the sets $T_0,\, T_1,\, T_2,\, \ldots$, where the set $T_n$ consists 
of all the reflections with $\infty$-height equal to $n$. 

These $T_n$'s were utilized to demonstrate 
nice regularity properties of $W$ (\cite[Ch. 5]{TE09}). Furthermore, they 
gave rise to a family of modules in the generic Iwahori-Hecke algebra associated to $W$, 
and in turn,  these modules  were  
used to prove a weak form of Lusztig's conjecture on the boundedness
of the $\mathbf{a}$-function (Dyer, unpublished). 
It is also known (Dyer, unpublished) that if $W$ is of finite rank, then there are finitely
many reflections in $T_n$ for each $n$.

In this section we prove that for an arbitrary reflection $t\in T$ whose $\infty$-height equals $n$, 
the corresponding positive root $\alpha_t$ dominates precisely $n$ other positive
roots. This observation will then establish a bijection between the set
of all reflections in $W$ with $\infty$-height equal to $n$ and the set of all positive roots
each dominates precisely $n$ other positive roots. Recent results on dominance obtained in \cite{FU1} 
may then be immediately applied to the  $T_n$'s,
answering a number of basic questions about these $T_n$'s. 
 
Following \cite{HT97} and \cite[$\S$~ 4.7]{ABFB}, we generalize the definition of dominance to the whole of $\Phi$
(whereas in \cite{BH93} and \cite{BB98}, dominance was only defined on $\Phi^+$), and we 
stress that all the notations are the same as in the previous section.

\begin{definition}
\label{def:dom}
\noindent\rm{(i)}\quad Let $W'$ be a reflection subgroup of $W$, and let $x, y\in
\Phi(W')$. Then we say that $x$ \emph{dominates} $y$ with respect to $W'$ if 
$$\{\,w\in W'\mid wx\in \Phi^-(W')\,\}\subseteq \{\,w\in W'\mid wy\in
\Phi^-(W')\,\}.$$ If $x$ dominates $y$ with respect to $W'$ then we write
$x\dom_{W'} y$. 

\noindent\rm{(ii)}\quad Let $W'$ be a reflection subgroup of $W$ and let $x\in
\Phi^+(W')$. Define
$D_{W'}(x) =\{\, y\in \Phi^+(W')\mid \text{$y\neq x$ and $x\dom_{W'} y$\,}\}$.
If $D_{W'}(x)=\emptyset$ then we call $x$ \emph{elementary with respect to}
$W'$. For each non-negative integer $n$, define $D_{W',n}=\{\,x\in \Phi^+(W')\mid
\#D_{W'}(x)=n\,\}$. In the case that $W'=W$, we write $D(x)$ and $D_n$ in place
of $D_{W'}(x)$ and $D_{W', n}$ respectively. If $D(x)=\emptyset$ then we call
$x$ \emph{elementary}. 
\end{definition}

It is readily checked that dominance with respect to any reflection subgroup $W'$ of
a Coxeter group $W$ is a partial ordering on $\Phi(W')$. 
The following lemma summarizes some basic properties of
dominance:
\begin{lemma}{\textup{(\cite[Lemma 2.2]{FU1})}}
  \label{lem:basicdom}
\rm{(i)}\quad Let $x, y \in \Phi^+$ be arbitrary. Then $x \dom_W y$ if and only
if $(x,y) \geq 1$ and $\dep(x)\geq \dep(y)$.

\noindent\rlap{\rm{(ii)}}\qquad Dominance is $W$-invariant, that is,  if $x \dom_W y$ then
 $wx \dom_W wy$ for all $w \in W$.

\noindent\rlap{\rm{(iii)}}\qquad Let $x, y \in \Phi$ be such that $x \dom_W y$.
Then
    $-y \dom_W -x$.

\noindent\rlap{\rm{(iv)}}\qquad Let $x,y\in \Phi$. Then there is dominance
between $x$ and $y$ if and only if $(x,y)\geq 1$. 
\qed
\end{lemma}

\begin{corollary}
 \label{co:basic}
Let $x, y\in \Phi$, and let $W'$ be an arbitrary reflection subgroup containing both
$r_x$ and $r_y$.

\noindent\rlap{\rm{(i)}}\qquad  There is dominance with respect to $W'$ between $x$ and $y$ if and
only if $(x, y)'\geq 1$, where $(\,,\,)'$ is the restriction of $(\,,\,)$ to the
subspace $\spa(\Delta(W'))$.

\noindent\rlap{\rm{(ii)}}\qquad 
$x\dom_W y$ if and only if $x\dom_{W'} y$.
\end{corollary}
\begin{proof}
 \rm{(i)}\quad  Follows from Lemma \ref{lem:basicdom} (iv) applied to the
Coxeter group $W'$ and the datum $\mathscr{C}'=(\,\spa(\Delta(W')),\,
\Delta(W'),\, (\,,\,)'\,)$.

\noindent\rlap{\rm{(ii)}}\qquad The desired result is trivially true if $x=y$, so 
we may assume that $x\neq y$. It is clear that $x\dom_W y$ implies that
$x\dom_{W'} y$. Conversely, suppose that $x\dom_{W'} y$. Then part (i) yields
that $(x,y)=(x, y)'\geq 1$. Thus Lemma~\ref{lem:basicdom}~(iv) yields that
either $x\dom_W y$, or else $y\dom_W x$. If the latter is the case, then by the
first part of the current proof, $y\dom_{W'} x$, and hence it follows that
$x=y$ (since dominance with respect to $W'$ is a partial ordering), 
contradicting our choice of $x$ and $y$. 
\end{proof}

Next is a well-known result whose proof can be found in the remarks
immediately before Lemma~2.3 of \cite{BH93}:
\begin{lemma}\textup{(\cite{BH93})}
 \label{lem:nodom}
There is no non-trivial dominance between positive roots in the root system of a finite Coxeter group.
\qed 
\end{lemma}

Then we have a technical result which is
going to be used repeatedly in the rest of this paper.
\begin{proposition}
 \label{pp:rootsys}
Let $\alpha, \beta\in \Phi^+$ with $(\alpha, \beta)\leq -1$, and 
let $W'$ be the dihedral reflection subgroup generated by $r_{\alpha}$ and $r_{\beta}$.
Further, if we set $\theta =\cosh^{-1}(-(\alpha, \beta))$, and for each $i\in \Z$ adopt the notation 
\begin{equation}
\label{eq:c_i}
c_i=
\begin{cases}
\frac{\sinh i\theta}{\sinh \theta},  \text{ if $\theta \neq 0$}\\
i, \text{ }\text{ }\text{ }\text{ }\text{ }\text{ }\,\text{if $\theta =0$}.
\end{cases}
\end{equation}
Then 

\noindent\rlap{\rm{(i)}}\qquad $W'$ is infinite, and 
$\Phi(W')=\{\,c_{i\pm 1} \alpha +c_i \beta\mid i\in \Z  \,\}$.

\noindent\rlap{\rm{(ii)}}\qquad Suppose that $x, y\in \Phi(W')$. Then $(x, y)\in (-\infty, -1]\cup [1, \infty)$, 
and in particular, if $x\neq \pm y$ then $\langle \{\,r_x,\,r_y\,\}\rangle$ is an infinite dihedral
reflection subgroup. More specifically, 
\begin{align*}
&\text{(a) If $x=c_{n+1}\alpha +c_n \beta$ and $y= c_{m+1}\alpha + c_m
\beta$, then either}\\
&\qquad\text{$(x,y)= \cosh((n-m)\theta) \geq 1$ if $\theta\neq 0$, or  $(x, y)=1$ if $\theta=0$.}\\
&\text{(b) If $x=c_{n+1}\alpha +c_n \beta$ and $y= c_{m-1}\alpha + c_m
\beta$, then either}\\
&\qquad\text{$(x,y)= -\cosh((n+m)\theta) \leq -1$ if $\theta\neq 0$, or $(x, y)=-1$ if $\theta=0$.}\\
&\text{(c) If $x=c_{n-1}\alpha +c_n \beta$ and $y= c_{m+1}\alpha + c_m
\beta$, then either}\\
&\qquad\text{$(x,y)= -\cosh((n+m)\theta) \leq -1$ if $\theta\neq 0$, or $(x, y)=-1$ if $\theta=0$.}\\
&\text{(d) If $x=c_{n-1}\alpha +c_n \beta$ and $y= c_{m-1}\alpha + c_m
\beta$, then either}\\
&\qquad\text{$(x,y)= \cosh((n-m)\theta) \geq 1$ if $\theta\neq 0$, or  $(x, y)=1$ if $\theta=0$.}
\end{align*}

\noindent\rlap{\rm{(iii)}}\qquad If $x\in \Phi^+(W')\setminus \{\,\alpha, \beta\,\}$ then $D_{W'}(x)\neq \emptyset$.
\end{proposition}

\begin{proof}
\rm{(i)} \quad  Proposition 4.5.4 (ii) of \cite{ABFB} implies that $W'$ is infinite, and the rest of
statement follows from direct calculations similar to those in Proposition~\ref{pp:anu1}.

\noindent\rlap{\rm{(ii)}}\qquad Follows from Part (i) above and a direct calculation.

\noindent\rlap{\rm{(iii)}}\qquad If $x\in \Phi^+(W')\setminus\{\,\alpha, \beta\,\}$ then Part (i) above yields that
either $x=c_{n+1}\alpha+c_n\beta$ (for some $n\neq 0$), or else $x=c_{n-1}\alpha+c_n\beta$ (for some $n\neq 1$). Then Part (ii) 
above and Corollary~\ref{co:basic}~(i) imply that we can find some $y\in \Phi^+(W')\setminus\{x\}$ such that $x\dom_{W'}y$. 

\end{proof}

The other key object to be studied in this section is the numeric function \emph{$\infty$-height} on $T$. 
As mentioned in the introduction section, this function is defined in terms of infinite dihedral reflection
subgroups of $W$, and in order to make a precise definition of this function we need a few technical results
on infinite dihedral reflection subgroups. We begin with the following well-known one, and for completeness, 
we include a proof here.

\begin{proposition}\textup{(Dyer \cite{bh})}
 \label{pp:dih}
Suppose that $\alpha, \beta\in \Phi^+$ are distinct. Let 
 $W'=\langle\,\{\,r_\gamma\mid \gamma\in (\R\alpha+\R\beta)\cap
\Phi^+\,\}\,\rangle$.
Then $W'$ is a dihedral reflection subgroup of $W$.
\end{proposition}
\begin{proof}
Suppose for a contradiction that $W'$ is not dihedral. Then $\#S(W')\geq 3$, and
let $x_1, x_2, x_3 \in \Delta(W')$ be distinct. Theorem~\ref{th:croots}~(ii)
then yields that $(x_i, x_j )\leq 0$ whenever $i,j\in\{1, 2, 3\}$ are different. 
Clearly $x_1, x_2, x_3$
are all in the two dimensional subspace $\R\alpha+\R\beta$, and thus a
contradiction arises if we could show that $x_1, x_2, x_3$ are linearly independent.
Let $c_1, c_2, c_3\in \R$ be such that $c_1 x_1+ c_2 x_2+c_3 x_3 = 0$. Since
$x_1, x_2, x_3 \in \Phi^+$, and $0\notin \PLC(\Delta)$,  it follows that $c_1, c_2,
c_3$ cannot be all positive or all negative. Rename $x_1, x_2, x_3$ if
necessary, we have the following three possibilities:
\begin{align}
 \label{eq:a}
&c_1, c_2 \geq 0 \quad\text{and}\quad c_3<0, \\
\noalign{\hbox{or}}
 \label{eq:b}
&c_1, c_2 \leq 0 \quad\text{and}\quad c_3>0,\\
\noalign{\hbox{or}}
\label{eq:c}
&c_1, c_2, c_3 =0. 
\end{align}
If (\ref{eq:a}) is the case then $0=(c_1x_1+c_2x_2+c_3x_3, x_3)<0$, and if
(\ref{eq:b}) is the case then $0=(c_1x_1+c_2x_2+c_3x_3, x_3)>0$, both are
clearly absurd. Hence (\ref{eq:c}) must be the case and $x_1, x_2, x_3$ are
linearly independent, a contradiction as required. 
\end{proof}

Let $\alpha, \beta\in \Phi^+$ be distinct. Let $W''$ be an arbitrary dihedral
reflection subgroup of $W$ containing the dihedral reflection subgroup $\langle
\{ r_\alpha, r_\beta\}\rangle$. Let $x, y$ be the canonical roots for $W''$. It
can be readily checked that $\R x+\R y=\R\alpha +\R\beta$, and hence $x, y\in
(\R\alpha+\R\beta)\cap \Phi^+$. It then follows that $W'' \subseteq
\langle\,\{\,r_\gamma \mid \gamma\in (\R\alpha+\R\beta)\cap \Phi^+\,\}\,\rangle$.
This observation together with Proposition~\ref{pp:dih} readily yield the
following well-known result:

\begin{proposition}
 \label{pp:maxdih}
Every dihedral reflection
subgroup $\langle \{r_\alpha, r_\beta\}\rangle$ of $W$ 
(where $\alpha, \beta\in \Phi^+$ are distinct),  
is contained in a unique maximal dihedral reflection subgroup, namely $\langle\,
\{\, r_\gamma \mid \gamma\in \Phi^+\cap(\R\alpha+\R\beta)\,\}\,\rangle$.
\qed
\end{proposition}

\begin{definition}
\label{def:h}

\rm{(i)}\, Define $\mathscr{M}$ to be the set of all maximal dihedral reflection
subgroups of $W$.

\noindent\rlap{\rm{(ii)}}\qquad Define $\mathscr{M}_{\infty}$ to be the set $\{\,W'\in
\mathscr{M}\mid \#W'=\infty\,\}$.

\noindent\rlap{\rm{(iii)}}\qquad For each $t\in T$, define
$\mathscr{M}_t$ to be the set $\{\,W'\in \mathscr{M}\mid t\in W'\,\}$. 

\noindent\rlap{\rm{(iv)}}\qquad Let $W'$ be a reflection subgroup of $W$, and let
$t\in W'\cap T$. Define the \emph{standard height}, $h_{(W',\,S(W')}(t)$, of $t$
with respect to the
Coxeter system $(W', S(W'))$ to be
$$\min\{\,\ell_{(W',\,S(W'))}(w)\mid \text{$w\in W'$,
$w\alpha_t\in \Delta(W')$}\,\}.$$
For the standard height of $t$ with respect to the Coxeter system $(W, S)$,  
we simply write $h(t)$ in place of $h_{(W, \, S)}(t)$. 
\end{definition}

\begin{remark}
For arbitrary reflection subgroup $W'$ of $W$, the depth function naturally applies
to $\Phi^+(W')$: if $x\in \Phi^+(W')$, then the \emph{depth} of $x$ relative to $S(W')$
(written $\dep_{(W',\, S(W'))}(x)$) is defined to be 
$$\min\{\,\ell_{(W', \, S(W'))}(w) \mid \text{$w\in W'$, and $wx\in \Phi^-(W')$}    \,\}.$$
Now for each $t\in W'\cap T$, it is easily checked that 
$$\dep_{(W',\, S(W'))}(\alpha_t)=h_{(W', \, S(W'))}(t)+1,$$
 and hence applying Lemma~\ref{lem:pre}~(i)
to the Coxeter system $(W', S(W'))$ yields that
\begin{equation}
\label{eq:height}
 h_{(W', \, S(W'))}(t)=\frac{\ell_{(W', \, S(W'))}(t)-1}{2}.
\end{equation}
\end{remark}

The following appears in \cite{TE09}, and for completeness we give a proof here:
\begin{lemma}
 \label{lem:t}
For each $t\in T$, we have $T\setminus\{t\}=\biguplus\limits_{W'\in \mathscr{M}_t}((W'\cap
T)\setminus\{t\})$.
\end{lemma}
\begin{proof}
 It is readily checked that $T\setminus\{t\}=\bigcup_{W'\in
\mathscr{M}_t}((W'\cap T)\setminus\{t\})$, and hence we only need to check that
this union is indeed disjoint. Suppose for a contradiction that there are distinct
$W_1, W_2\in \mathscr{M}_t$ with $r\in W_1\cap W_2$ for some $r\in
T\setminus\{t\}$. Then clearly $\langle \{r, t\}\rangle \subseteq W_1$ and
$\langle \{r, t\}\rangle \subseteq~W_2$, contradicting
Proposition~\ref{pp:maxdih}.
\end{proof}

The canonical bijection $\psi\colon T \leftrightarrow \Phi^+$ and the above
immediately yield that:
\begin{corollary}
\label{co:root}
$\Phi^+\setminus \{\alpha\}=\biguplus\limits_{W'\in
\mathscr{M}_{r_{\alpha}}}(\Phi^+(W')\setminus \{\alpha\})$, for each $\alpha~\in~\Phi^+$. 
\qed
\end{corollary}

\begin{remark}
 \label{rm:l10}
In particular, the above corollary yields implies that for $t\in T$, if 
$W_1, W_2\in \mathscr{M}_t$ are distinct then $\Phi^+(W_1)\cap\Phi^+(W_2)=\{\alpha_t\}$.
\end{remark}

\begin{lemma}\textup{(\cite{TE09})}
 \label{lem:height}
Let $t\in T$ be arbitrary. Then 
$$h(t)=\sum\limits_{W'\in \mathscr{M}_t}h_{(W', \ S(W'))}(t).$$
\end{lemma}
\begin{proof}
For any reflection $t\in T$, Corollary \ref{co:root} yields that
\begin{equation}
\label{eq:e1}
\{\alpha\in \Phi^+\mid t\alpha\in
\Phi^-\}=\{\alpha_t\}\cup(\biguplus\limits_{W'\in \mathscr{M}_t}\{\alpha\in
\Phi^+(W')\setminus\{\alpha_t\}\mid t\alpha\in \Phi^-\,\}).
\end{equation} 
Since $h(t)=\frac{1}{2}(\ell(t)-1)=\frac{1}{2} (\#N(t)-1)$, it follows from
(\ref{eq:e1}) that
 \begin{align*}
 h(t)&=\frac{1}{2}(\sum\limits_{W'\in \mathscr{M}_t}\#\{\alpha \in
\Phi^+(W')\setminus\{\alpha_t\}\mid t\alpha \in \Phi^-(W')\,\})\\
     &=\sum\limits_{W'\in
\mathscr{M}_t}\frac{1}{2}(\ell_{(W',\,S(W'))}(t)-1)\qquad\text{( by (\ref{eq:ell}) ) }\\
     &=\sum\limits_{W'\in \mathscr{M}_t}
h_{(W',\,S(W'))}(t)\qquad\qquad\quad\text{( by (\ref{eq:height}) ) }.
\end{align*}
\end{proof}

\begin{definition}\textup{(\cite{TE09})}
 For $t\in T$, define the \emph{$\infty$-height} of $t$ to be 
$$h^{\infty}(t)=\sum\limits_{W'\in\mathscr{M}_t\cap \mathscr{M}_{\infty}}
h_{(W', \ S(W'))}(t), $$
and for each non-negative integer $n$, we define
$$T_n=\{t\in T\mid h^{\infty}(t)=n\}.$$
\end{definition}

Observe that from the above definition, it is not clear whether, for a specific
non-negative integer $n$, there is any reflection $t\in T$ with $h^{\infty}(t)=n$.
It turns out that a number of basic questions like this can in fact be resolved with the 
aid of the results obtained in \cite{FU1} once we prove the following: 
\begin{theorem}
 \label{th:bij}
For each non-negative integer $n$, there is a bijection $T_n\leftrightarrow D_n$
given by $t\leftrightarrow \alpha_t$. 
\end{theorem}
The proof of the above theorem will be deferred until we have all the necessary tools.

\begin{proposition}
 \label{pp:p6}
Suppose that $t\in T$, and let $W'$ be an infinite dihedral reflection subgroup containing
$t$. If $h_{(W', \ S(W'))}(t)\geq 1$ then there exists some $x\in \Phi^+(W')$
with $\alpha_t\dom_W x$.
\end{proposition}
\begin{proof}
Observe that the condition $h_{(W', \ S(W'))}(t)\geq 1$ is equivalent to 
$\alpha_t\notin \Delta(W')$, and hence the required result follows immediately 
from Proposition \ref{pp:rootsys}~(iii).
\end{proof}

The following proposition will be a key step to prove Theorem \ref{th:bij}:

\begin{proposition}
 \label{pp:chain}
Let $W'$ be an infinite dihedral reflection subgroup, and let $\Delta(W')=\{\,\alpha, \beta\,\}$.

\noindent\rm{(i)}\quad There are two disjoint dominance chains in $\Phi(W')$, namely:
\begin{multline}
\label{eq:chain1}
\cdots \dom_W \ r_\alpha r_\beta r_\alpha (\beta) \ \dom_W \ r_\alpha r_\beta
(\alpha) \ \dom_W r_\alpha (\beta) \ \dom_W  \ \alpha \\
 \dom_W \ (-\beta) \ \dom_W \ r_\beta(-\alpha) \  \dom_W \ r_\beta r_\alpha
(-\beta) \ \dom_W  \ \cdots 
\end{multline}
 and 
\begin{multline}
\label{eq:chain2}
\cdots  \dom_W \ r_\beta r_\alpha r_\beta (\alpha) \ \dom_W \ r_\beta r_\alpha
(\beta) \ \dom_W \ r_\beta (\alpha) \ \dom_W \ \beta \\
\dom_W (-\alpha) \ \dom_W \ r_\alpha (-\beta) \ \dom_W \ r_\alpha r_\beta
(-\alpha) \ \dom_W \ \cdots.
\end{multline}
In particular, each root in $\Phi(W')$ lies in exactly one of the above two chains, and the negative of 
any element of one chain lies in the other. Furthermore, the roots in $\Phi(W')$ dominated by either 
$\alpha$ or $\beta$ are all negative.

\noindent\rm{(ii)}\quad If $x\in \Phi(W')$ then $\#D_{W'}(x)=h_{(W',\, S(W'))}(r_x)$.
\end{proposition}
\begin{proof}
 \rm{(i)}\quad Theorem \ref{th:croots}~(ii) and \cite[Proposition 4.5.4 (ii)]{ABFB} yield that 
$(\alpha, \beta)\leq -1$. Hence it follows from Lemma~\ref{lem:basicdom}~(iv) that 
$\alpha \dom_W -\beta$ and $\beta\dom_W -\alpha$. Then we can immediately verify the existence of 
the two dominance chains (\ref{eq:chain1}) and (\ref{eq:chain2}), and from these two chains 
the remaining statements in part (i) follow readily.

\noindent\rm{(ii)}\quad Follows immediately from the definition of $h_{(W',\, S(W'))}(r_x)$
and the two dominance chains (\ref{eq:chain1}) and (\ref{eq:chain2}).
\end{proof}

\begin{proposition}
 \label{pp:p7}
Suppose that $x, y\in \Phi^+$ are distinct with $x\dom_W y$, and let $W'$ be a
dihedral reflection subgroup containing $r_x$ and $r_y$. Then $h_{(W', \ S(W'))}(r_x)\geq 1$.

\end{proposition}
\begin{proof}
It follows from Corollary \ref{co:basic} (ii) that $x\dom_{W'} y$, so Lemma
\ref{lem:nodom} above yields that $W'$ is an infinite dihedral reflection subgroup. Let
$\{\alpha, \beta\}=\Delta(W')$. We know from Proposition~\ref{pp:chain}~(i) that the roots
in $\Phi(W')$ dominated by either $\alpha$ or $\beta$ are all negative, and since 
$x\dom_W y\in \Phi^+$, it follows that $x\notin\{\,\alpha,\, \beta\,\}$. Hence by definition
$h_{(W', \ S(W'))}(r_x)\geq 1$.

\end{proof}

From the last two propositions we may deduce the following special case of
Theorem \ref{th:bij}:
\begin{lemma}
 There is a bijection $T_0\leftrightarrow D_0$ given by $t\leftrightarrow
\alpha_t$.
\end{lemma}
\begin{proof}
Let $t\in T_0$, and suppose for a contradiction that $\alpha_t\notin D_0$. Then
there exists $s\in T\setminus\{t\}$ such that $\alpha_t\dom_W \alpha_s$. Let
$W'$ be the unique maximal dihedral reflection subgroup of $W$ containing
$\langle \{s, t\}\rangle$. Proposition~\ref{pp:p7} yields that $h_{(W', \
S(W'))}(t)\geq 1$. Since $\alpha_t\dom_W \alpha_s$, it follows from 
Lemma~\ref{lem:nodom} that $W'\in \mathscr{M}_{\infty}$, and
consequently $h^{\infty}(t)\geq 1$, contradicting the
assumption that $t\in T_0$.

Conversely, suppose that $\alpha_t\in D_0$, and suppose for a contradiction that
$t\notin T_0$. Then there exists some $W'\in
\mathscr{M}_t\cap\mathscr{M}_{\infty}$ with $h_{(W', \ S(W'))}(t)\geq 1$. But
then Proposition~\ref{pp:p6} yields that $\alpha_t\notin D_0$, producing a
contradiction as required.
\end{proof}

Observe that Proposition \ref{pp:chain} (ii) can be equivalently stated as:
\begin{proposition}
 \label{pp:p11}
Suppose that $t\in T$, and suppose that $W'$ is an infinite dihedral reflection subgroup 
containing $t$. Then 
$$\#D_{W'}(\alpha_t)=h_{(W', \ S(W'))}(t).$$

\qed
\end{proposition}

\begin{proposition}
 \label{pp:p12}
Suppose that $t\in T$ is arbitrary. Then 
$$\biguplus\limits_{W'\in
\mathscr{M}_t\cap\mathscr{M}_{\infty}}D_{W'}(\alpha_t)=D(\alpha_t).$$
\end{proposition}
\begin{proof}
First we observe that Remark \ref{rm:l10} yields that the union of the sets
$D_{W'}(\alpha_t)$ over all $W'$ in $\mathscr{M}_t\cap\mathscr{M}_\infty$ is
indeed disjoint.

It is clear that 
$\biguplus\limits_{W'\in
\mathscr{M}_t\cap\mathscr{M}_{\infty}}D_{W'}(\alpha_t)\subseteq D(\alpha_t)$.

Conversely, suppose that $x\in D(\alpha_t)$. Let $W'$ be the unique maximal
dihedral reflection subgroup of $W$ containing $\langle \{ t, r_x\}\rangle$.
Then Corollary~\ref{co:basic}~(ii) yields that $\alpha_t\dom_{W'} x$.  Finally
since there is no non-trivial dominance in any finite Coxeter group, it follows
that $W'\in \mathscr{M}_\infty$, as required.
\end{proof}

Now we prove that for any reflection  $t\in W$, its $\infty$-height $h^{\infty}(t)$ 
equals the the number of positive roots strictly dominated by $\alpha_t$:  
\begin{theorem}
 \label{th:equal}
Let $t\in T$ be arbitrary. Then $h^{\infty}(t)=\#D(\alpha_t)$.
\end{theorem}
\begin{proof}
 It follows from Proposition \ref{pp:p11} and Proposition~\ref{pp:p12} that 
$$
 h^{\infty}(t)=\sum\limits_{W'\in
\mathscr{M}_{t}\cap\mathscr{M}_{\infty}}h_{(W', \ S(W'))}(t)= \sum\limits_{W'\in
\mathscr{M}_{t}\cap\mathscr{M}_{\infty}}\#D_{W'}(\alpha_t)
=\#D(\alpha_t).
$$
\end{proof}

Finally we are in a position to prove Theorem \ref{th:bij}:
\begin{proof}[Proof of Theorem \ref{th:bij}]
 The desired result follows immediately from Theorem \ref{th:equal}.
\end{proof}

Now combining Theorem 3.8 of \cite{FU1},  Corollary 3.9 of \cite{FU1},
Corollary 3.21 of \cite{FU1} and Theorem \ref{th:bij} above we may deduce:
\begin{corollary}
\label{co:T_n}
 \begin{itemize}
 \item [(i)] For each positive integer $n$, 
$$
T_n\subseteq \{\,tt't\mid\text{$t\in T_0$ and $t'\in T_m$ for some $m\leq
n-1$}\,\}.
$$ 
\item [(ii)] Suppose that $W$ is an infinite Coxeter group with $\# S<\infty$. 
Then $0<\#T_n\leq \ (\#T_0)^{n+1}-(\#T_0)^n$ 
for each positive integer $n$.
\end{itemize}
\qed
\end{corollary}

\begin{remark}
An upper bound for $\#T_0(=\#D_0)$ is given in \cite{BH93}, furthermore, for 
any fixed finitely generated Coxeter group, this number can be explicitly calculated 
following the methods presented in \cite{BB98}.
\end{remark}

\section{Dominance and Imaginary Cone}
\label{sec:IMC}

Kac introduced the the concept of an \emph{imaginary cone} in the study of the
imaginary roots of Kac-Moody Lie algebras. In \cite[Ch. 5]{VK} the imaginary
cone of a Kac-Moody Lie algebra was defined to be the positive cone 
on the positive imaginary roots. The 
generalization of imaginary cones to arbitrary Coxeter groups was 
first introduced by H\'ee in \cite{HE1}, and subsequently reproduced in \cite{HE2}.  
This generalization has also been studied by Dyer (\cite{MD12}) and Edgar (\cite{TE09}).
In this section we investigate the connections between this generalized 
imaginary cone and dominance in Coxeter groups, in particular, we show that
whenever $x$ and $y$ are roots of a Coxeter group, then $x \dom_W y$ if and 
only if  $x-y$ lies in the imaginary cone of that Coxeter group.

Let $(W, S)$ be the abstract Coxeter system associated to the Coxeter datum
$\mathscr{C}=(\,V,
\Delta, (\,,\,)\,)$ and let $\Phi$ be the corresponding root system. For any real
vector space $X$ we write $X^*=\Hom(X, \R)$. In this section
we take $X$ to be some suitable subspace of $V$. Also in this paper all \emph{cones}
are assumed to be convex cones. 
For any cone $C$ in $X$, we define 
$C^* = \{\, f\in X^* \mid f(v) \geq 0 \text{ for all } v\in C \,\}$ 
and call it the \emph{dual} of $C$; and for any cone $F$
in $X^*$, we define 
$F^*=\{\, v \in X \mid f(v)\geq 0 \text{ for all } f \in F\,\}$ 
and call it the \emph{dual} of $F$. If $W$ acts on $X$, then $X^*$ 
bears the \emph{contragredient representation} of $W$ in the following way:
if $w\in W$ and $f\in X^*$ then $wf\in X^*$ is given by the rule
$(wf)(v)=f(w^{-1}v)$ for all $v\in X$. It is readily checked that for a
cone $C$ in $X$ we have $C\subseteq C^{**}$, and also for any $w\in W$, we
have  $(w C)^*=w C^*$.

The following is a well-known result whose proof can be found in 
\cite[Notes (c), Lecture 1]{RB96}:
\begin{lemma}
 \label{lem:basic}
Suppose that $X$ is a real vector space of finite dimension, and let $C$ be a
cone in $X$. Then $(C^*)^*=\overline{C}$, where $\overline{C}$ is the
topological closure of $C$ in $X$ (with respect to the standard topology on
$X$).
\qed
\end{lemma}

Set $P= \PLC(\Delta) \cup \{ 0 \}$. It is clear that $P$ is a cone in $V$.
We define the Tits cone of  $W$ in the same way as in 5.13 of \cite{HM}:
\begin{definition}
The \emph{Tits cone} of the Coxeter group $W$
is the $W$-invariant set $U = \bigcup_{w \in W} w P^*$. 
\end{definition}

It is not obvious from the above definition that the Tits cone is indeed a cone, 
however, this can be made clear by the following result:
\begin{proposition}
\begin{equation}
\label{eq:Tits}
U=\{\, f\in \spa(\Delta)^* \mid \text{$f(x)\geq 0$ for all but finitely 
many $x\in \Phi^+$}     \,\}.
\end{equation}
\end{proposition}
\begin{proof}
Denote the set on the right hand side of (\ref{eq:Tits}) by $Y$, and for 
each $f\in \spa(\Delta)^*$ define $\n(f)$ by 
$\n(f)=\{\,x\in \Phi^+\mid f(x)<0\,\}$.

If $f\in U$ then $f=wg$ for some $w\in W$ and $g\in P^*$, and it is readily 
checked that $\n(f)\subseteq N(w^{-1})$. Since $N(w^{-1})$ is a finite set, 
it follows that $f\in Y$, and hence $U\subseteq Y$.
Conversely, suppose that $f\in Y$. If $\n(f)=\emptyset$ 
then $f\in P^*\subseteq U$. Thus we may assume that 
$\#\n(f)>0$, and proceed with an induction. Observe that
then there exists some $\alpha \in \Delta$ such that 
$f(\alpha)<0$. It is then readily checked that 
$\#\n(r_{\alpha} f)=\#\n(f)-1$, and hence it follows from the 
inductive hypothesis that $r_{\alpha}f\in U$. Since $U$ is 
$W$-invariant, it follows that $f\in U$, and hence $Y\subseteq U$.
\end{proof}

\begin{lemma}
\label{lem:tits}
$U^*=\bigcap\limits_{w\in W} w (P^*)^*$.
Furthermore,
$U^* = \bigcap\limits_{w \in W} w P$, whenever $\Delta$ is a finite set.
\end{lemma}

\begin{proof}
\begin{align}
 \label{eq:tits}
U^* &= \{\, v\in V \mid f(v)\geq 0 \text{, for all $f \in U$}\,\}\nonumber\\  
&=  \{\, v\in V \mid (w\phi)(v) \geq 0, \  \text{for all $\phi \in P^*$, and
for all $w\in W$}     \,\}\nonumber \\
&=\{\,v\in V \mid \phi(w^{-1}v)\geq 0, \  \text{for all $\phi \in P^*$, and for
all $w\in W$}  \,\} \nonumber\\
&=\bigcap_{w \in W} \{\, v \in V \mid \phi(w^{-1}v)\geq 0, \ \text{for all
$\phi \in P^*$} \,\}\nonumber\\
&=\bigcap_{w \in W} \{\,wv \in V \mid \phi(v) \geq 0, \ \text{for all $\phi
\in P^*$}   \,\}\nonumber\\
&= \bigcap_{w \in W} \{\,wv \in V \mid v \in  (P^*)^*  \,\}.
\end{align}
Let $X=\spa(\Delta)$. If $\#\Delta$ is finite then it follows from Lemma~\ref{lem:basic}
 that $(P^*)^*=\overline{P}$. It is clear that $P$ is topologically closed, hence
(\ref{eq:tits}) yields that $U^*=\bigcap_{w\in W} wP$ when $\Delta$ is a finite set.
\end{proof}

\begin{lemma}
\label{lem:l5}
Suppose that $v \in V$ has the property that $(a, v) \leq 0$ for all $a\in \Delta$.
Then $wv-v \in P$ for all $w \in W$. Moreover, if $v\in P$ then $v \in U^*$.
\end{lemma}
\begin{proof}
Use induction on $\ell(w)$. Note that if $\ell(w)=0$ then there is nothing to
prove. If $\ell(w) \geq 1$ then we may write $w=w' r_a$ where $w'\in W$ and $a
\in \Delta$ with $\ell(w)= \ell(w')+1$. Then Proposition~\ref{pp:anu2} yields that 
$w'a\in \Phi^+\subseteq P$, and we have
\begin{displaymath}
\begin{split} 
wv-v = (w'r_a) v-v &= w'(v-2(v,a)a)-v \\
                   & = (w'v-v) -2(a,v) w'a.
\end{split}
\end{displaymath}
Note that by the inductive hypothesis $w'v-v \in P$. Since
$(a, v)\leq 0$, it follows from the above that $wv-v\in P$.

If $ v \in P$ then $wv = (wv-v) +v \in P$ for all $w \in W$, and hence $v \in
\bigcap\limits_{w \in W} w^{-1} P\subseteq U^*$.
\end{proof}

The following is a useful result from \cite{FU1}:
\begin{proposition}\textup{(\cite[Proposition 3.4]{FU1})}
 \label{pp:candih}
Suppose that $x, y\in \Phi$ are distinct with $x\dom_W y$. Let $W'$ be the dihedral
reflection subgroup generated by $r_x$ and $r_y$, and let 
$\Delta(W')= \{\,\alpha, \beta\,\}$. Then there exists some $w\in W'$ such that either
\begin{equation*}
\left\{
\begin{array}{rl}
w x &= \alpha\\
w y &= -\beta
\end{array} \right.\qquad\text{or else}\qquad
\left\{
\begin{array}{rl}
w x&= \beta\\
w y&= -\alpha.
\end{array} \right.
\end{equation*}
In particular, $(x, y)=-(a, b)$.
\qed
\end{proposition}

\begin{proposition}
\label{pp:tits}
Suppose that $x,y \in \Phi$ such that $x \dom_W y$. Then  $w(x-y) \in
\PLC(\Delta)$ for all $w \in
W$, that is, $x-y \in U^*$.
\end{proposition}

\begin{proof}
The assertion is trivially true if $x=y$, so we may assume that $x\neq y$.
Since $x\dom_W y$, Lemma~\ref{lem:basicdom}~(iv) yields that $(x, y)\geq 1$.
Let $W'$ be the (infinite) dihedral subgroup of $W$ generated by $r_x$ and
$r_y$. Let $S(W') = \{ s,t \}$ and $\bigtriangleup(W') = \{ \alpha_s, \alpha_t
\}$. Proposition~\ref{pp:candih} yields that $(\alpha_s, \alpha_t)=-(x,
y)\leq -1$. 
Set $c_i$ as in Proposition~\ref{pp:rootsys}
for each $i\in \Z$.
Since $x\dom_W y$, it follows that $(x, y)\geq 1$, and 
Proposition~\ref{pp:rootsys}~(ii) then yields that either
\begin{equation*}
 \left\{
\begin{array}{rl}
x &= c_{n+1}\alpha_s + c_n\alpha_t\\
y &= c_{m+1}\alpha_s + c_m \alpha_t
\end{array} \right. \qquad\text{or else}\qquad
\left\{
\begin{array}{rl}
x &= c_{n-1}\alpha_s + c_n\alpha_t\\
y &= c_{m-1}\alpha_s + c_m \alpha_t
\end{array} \right..
\end{equation*}
Next we shall show that $n > m$. Suppose for a contradiction that $m\geq n$.
Then either $x =y$ (when $n=m$) or else there will be a $w\in W'$ such that $w x
\in \Phi(W') \cap \Phi^-$ and yet $wy \in \Phi(W') \cap \Phi^+$ (when $n<m$),
both contradicting the fact that $x\dom_W y$. Since $c_n>c_m$ whenever $n>m$, it
follows that $x-y \in \PLC(\Delta)$. Given the $W$-invariance of dominance, 
for any $w \in W$, repeat the above argument with $x$ replaced by $wx$ and $y$
replaced by $wy$, we may conclude that $w(x-y) \in \PLC(\Delta)\subseteq (P^*)^*$.
It then follows from Lemma~\ref{lem:tits} that $x-y\in U^*$.
\end{proof}
 
When $\#\Delta$ is finite, it can be checked that Lemma ~\ref{lem:tits}
yields that whenever $x, y\in \Phi$ such that $x-y\in U^*$, then $x\dom_W y$.
In fact we can remove this finiteness condition and still prove the same result, and 
to do so we need  some special notations and few extra elementary results. We thank the 
referee of this paper for prompting us to look into this direction.

\begin{notations}
For a subset $I$ of $S$ we set
$\Delta _I=\{\, x\in \Delta \mid r_x\in I\,\}$;
$V_I=\spa (\Delta_I)$;
$W_I= \langle\, I\,\rangle$; and
$P_I=\PLC(\Delta_I)\cup \{0\}$. Furthermore, we set
\begin{align*}
P_I^*&=\{\, f\in \Hom(V_I, \R)\mid \text{$f(x)\geq 0$ for all $x\in P_I$}\,\};\\
\noalign{\hbox{and}}
P_I^{**}&=\{\, x\in V_I \mid \text{$f(x)\geq 0$ for all $f\in P_I^*$}\,\}.
\end{align*}
\end{notations}
Then $\mathscr{C}_I=(\, V_I, \Delta_I, (\,,\,)_I\,)$ (where $(\,,\,)_I$ is the 
restriction of $(\,,\,)$ on $V_I$) is a Coxeter datum 
with corresponding Coxeter system $(W_I,I)$, and we call $W_I$ the 
\emph{standard parabolic subgroup} of $W$ corresponding to $I$. Clearly
$W_I$ preserves $V_I$.

\begin{lemma}
 \label{lem:P_I}
Suppose that $I$ is a subset of $S$. Then $P^{**}\cap V_I\subseteq P_I^{**}$.
\end{lemma}
\begin{proof}
 Write $V=V_I\oplus V_I'$, where $V_I'$ is a vector space complement
of $V_I$. Consequently, every $v\in V$ is uniquely written as 
$v=v_I+v_I'$, where $v_I\in V_I$ and $v_I'\in V_I'$. Then we observe that 
every $g\in P_I^*$ gives rise to a $g'\in P^*$ as follows: for any $v\in V$, 
simply set $g'(v)=g(v_I)$. Now let $x\in P^{**}\cap V_I$ and $f\in P_I^*$ be 
arbitrary. Then $f(x)=f'(x)\geq 0$, since $f'\in P^*$ and $x\in P^{**}$. Hence 
$x\in P_I^{**}$, and so $P^{**}\cap V_I\subseteq P_I^{**}$.
\end{proof}

\begin{proposition}
 \label{pp:TitsConverse}
Let $x, y\in \Phi$. Then $x-y\in U^*$ if and only if $x\dom_W y$.
\end{proposition}
\begin{proof}
 By Proposition \ref{pp:tits} we only need to prove that when $x$ and $y$ are both roots then
 $x-y\in U^*$ implies that $x\dom_W y$. The assertion certainly 
holds if $x=y$, thus we only need to check the case
when $x\neq y$.

Since dominance and $U^*$ are both $W$-invariant, it follows that we only need 
to prove the following statement: if $x\in \Phi^-$ then $y\in \Phi^-$ too. 

Take $I=\{\, r_{\alpha}\mid \alpha \in \supp(x)\cup \supp(y)\,\}$, and note 
that in particular, $I$ is a finite set. Now in view of Lemma~\ref{lem:tits},
Lemma~\ref{lem:P_I} and the fact that $W_I$ preserves $V_I$ we have
\begin{align*}
 x-y\in (\bigcap_{w\in W} wP^{**})\cap V_I\subseteq 
      (\bigcap_{w\in W_I} wP^{**})\cap V_I &\subseteq
       \bigcap_{w\in W_I} w (P^{**}\cap V_I)\\
  &\subseteq \bigcap_{w\in W_I} wP_I^{**}\\
  &= \bigcap_{w\in W_I} w P_I,
\end{align*}
where the equality follows from Lemma \ref{lem:basic}, since $I$ is a finite set.
Thus $x-y\in P_I$, and this implies, precisely, that $y\in \Phi^-$ whenever $x\in \Phi^-$.
\end{proof}

Next we have a technical result which is a key component of the main 
theorem of this section.

\begin{proposition}
 \label{pp:key}
Suppose that $x, y\in \Phi$ are distinct with $x\dom_W y$. Then there exists
some $w\in W$ such that $wx\in \Phi^+$, $wy\in \Phi^-$ and $(w(x-y), z)\leq 0$
for all $z\in \Phi^+$.
\end{proposition}
\begin{proof}
Clearly it is enough to show that under such assumptions there exists some
$w\in W$ with $wx\in \Phi^+$, $wy\in \Phi^-$ and $(w(x-y), z)\leq 0$ for all
$z\in \Delta$. 

Let $W'$ be the (infinite) dihedral reflection subgroup of $W$ generated by $r_x$ and
$r_y$, and let $\bigtriangleup(W') =\{ a_0, b_0 \}$.  Clearly $a_0$, $b_0 \in\Phi^+$,
and Proposition~\ref{pp:candih} yields that $(a_0, b_0) =-(x,y) \leq-1$, 
furthermore, there is some $u \in \langle \{r_x, r_y\} \rangle$ such that
either 
\begin{equation}
\label{eq:pair}
\left\{
\begin{array}{rl}
u(x)&= a_0\\
u(y)&= -b_0
\end{array} \right.\qquad\text{or else}\qquad
\left\{
\begin{array}{rl}
u(x)&= b_0\\
u(y)&= -a_0.
\end{array} \right.
\end{equation}
At any rate, $u(x-y) = a_0 + b_0$. Since the $W$-action preserves $(\,,\,)$, it follows that 
$(a_0, a_0)=1=(b_0, b_0)$, and hence
$(a_0+b_0, a_0) \leq 0$ and $(a_0 + b_0 , b_0 )\leq 0$. However
there may exist some $c_1 \in \Delta$ with $(a_0+b_0, c_1 ) > 0$. If this is the
case, set $a_1 = r_{c_1} a_0$ and $b_1 = r_{c_1} b_0$.
Recall that $(d, c_1) \leq 0$ for all $d\in \Delta\setminus\{c_1\}$, so it
follows that
\begin{equation}
\label{eq:supp}
c_1\in \supp(a_0)\cup \supp(b_0).
\end{equation}
Since $(a_0+b_0 , c_1) >0$, whereas $(a_0+b_0, a_0) \leq 0$ and $(a_0+b_0, b_0)
\leq 0$, it follows that $a_0\neq c_1$ and $b_0 \neq c_1$. Therefore we see that
$a_1, b_1 \in \Phi^+$, and $(a_1, b_1) =(a_0, b_0)\leq -1$.  Consequently
Theorem~\ref{th:croots}~(ii) yields that $a_1, b_1$ are the canonical roots for
the root subsystem $\Phi(\langle \{r_{a_1}, r_{b_1}\}\rangle)$. 
Since $r_{c_1}(a_0+b_0)=a_0+b_0-2(a_0+b_0, c_1)c_1$ and $(a_0+b_0, c_1)>0$, it
follows  that
$$
\supp(a_1)\cup\supp(b_1) \subseteq \supp(a_0)\cup \supp(b_0),
$$ 
and
$$\sum\limits_{a\in \Delta} \coeff_a (a_1) + \sum\limits_{a\in \Delta} \coeff_a (b_1)
<\sum\limits_{a\in \Delta} \coeff_a (a_0) + \sum\limits_{a\in \Delta} \coeff_a
(b_0).$$ Moreover, since $(a_0+b_0, c_1 ) >0$, it follows that at least one of
$(a_0, c_1)$ or $(b_0, c_1)$ must be strictly positive. Hence 
Lemma~\ref{lem:pre} yields that
$$dp(a_1) + dp(b_1) \leq dp(a_0)+dp(b_0).$$

Repeat this process and we can obtain new pairs of positive roots $\{\, a_2, b_2
\,\}, \ldots, \{\,a_{m-1}, b_{m-1}\,\}, \{\,a_m, b_m\,\}$ with 
\begin{align*}
\supp(a_m)\cup \supp(b_m) &\subseteq \supp(a_{m-1})\cup\supp(b_{m-1})\subseteq
\cdots \\
&\subseteq \supp(a_0)\cup\supp(b_0)
\end{align*}
and $\dep(a_m)+\dep(b_m)\leq \dep(a_{m-1})+\dep(b_{m-1})\leq \cdots \leq
\dep(a_0)+\dep(b_0)$, 
so long as we can find a $c_m \in \Delta$ such that $(a_{m-1}+b_{m-1}, c_m ) >0$.
Note that this process only terminates at a pair $\{\, a_n, b_n\,\}$ for some
$n$, if $(a_n+b_n, z) \leq 0$ for all $z \in \Delta$. Now if we could
show that this process terminates at some such $\{a_n, b_n\}$ after a finite
number of iterations, then we have in fact
found a $w\in W$ given by 
\begin{equation}
 \label{eq:w}
w~=~r_{c_n} r_{c_{n-1}} \cdots r_{c_1} u, \text{   where $u$ is as in (\ref{eq:pair})}, 
\end{equation}
satisfying 
$$(w(x-y), z) = (r_{c_n}\cdots r_{c_1} (a_0+b_0), z)=(a_n+b_n, z)\leq 0$$ for
all $z\in \Delta$.

Observe that the set of positive roots having depth less than the specific bound 
$\dep(a_0)+\dep(b_0)$ and support in a fixed finite subset 
$\supp(a_0)\cup\supp(b_0)$ of $\Delta$ is finite, indeed, Lemma~\ref{lem:pre}~(ii) 
implies that there are at most $\#(supp(a_0)\cup\supp(b_0))^{\dep(a_0)+\dep(b_0)}$
many such positive roots.  Hence it  follows that the
possible pairs of positive roots $\{a_i, b_i\}$ obtainable in the above process must
be
finite too. Finally since  
$$\sum\limits_{a\in \Delta} \coeff_a (a_j) + \sum\limits_{a\in \Delta} \coeff_a (b_j)
<\sum\limits_{a\in \Delta} \coeff_a (a_i) + \sum\limits_{a\in \Delta} \coeff_a
(b_i)$$ 
for all $j>i$, it follows that the sequence $\{ a_0, b_0\}, \{a_1, b_1\}, \cdots $
must terminate at $\{a_n, b_n\}$ for some finite $n$, as required.

Finally, keep $w$ as in (\ref{eq:w}), we see from the above construction that either $wx=a_n\in \Phi^+$
and $w y=-b_n\in \Phi^-$, or else $w x= b_n\in \Phi^+$ and $w y=-a_n\in \Phi^-$.

\end{proof}

\begin{definition}
We define the \emph{imaginary cone} $Q$ of $W$ by 
$$
Q =\{\, v\in U^* \mid (v,a) \leq 0 \text{ for all but finitely many } a\in
\Phi^+ \,\}.
$$
\end{definition}

The following result was obtained independently by Dyer as a 
consequence of \cite[Theorem 6.3]{MD12}, stating that \emph{the imaginary cone of a 
reflection subgroup is contained in that of the over-group}. 
\begin{theorem}
\label{th:cone}
\label{th:imc}
Suppose that $x$, $y \in \Phi$ such that $x \dom_W y$. Then $x-y\in Q$. 
\end{theorem}

\begin{proof}
By Proposition \ref{pp:tits} we know that $x-y\in U^*$, thus to prove the
desired
result, we only need to show that $(x-y, z)\leq 0$ for all but finitely many
$z\in \Phi^+$. Suppose that $z\in \Phi^+$ such that $(x-y, z)>0$. Let $w\in W$
be as in Proposition~\ref{pp:key}. Then $(w(x-y), wz)>0$, and by
Proposition~\ref{pp:key} this is possible only if $z\in N(w)$. Since
$\#N(w)$ is clearly finite (of size $\ell(w)$), it follows that indeed
$(x-y, z)\leq 0$ for
all but finitely many $z\in \Phi^+$.  
\end{proof}

\begin{remark}
The above theorem is a special case of Dyer's result when the subgroup is dihedral.
In fact, Dyer's result, when applied to dihedral reflection subgroups, implies that 
if $x$ and $y$ are roots with $x\dom_W y$ then $x-cy\in Q$ for an explicit range of 
$c\in \R$ depending on the value of $(x, y)$. Our formulation was first suggested to 
us by Howlett and Dyer, and we gratefully acknowledge their help. 
\end{remark}

Theorem \ref{th:imc} combined with Proposition \ref{pp:TitsConverse}
immediately imply the following:

\begin{corollary}
 \label{co:imc}
Let $x, y\in \Phi$. Then $x-y\in Q$ if and only if $x\dom_W y$.
\qed
\end{corollary}

\begin{remark}
Incidentally, we observe from Proposition \ref{pp:TitsConverse} and 
Corollary~\ref{co:imc} that when $x, y\in \Phi$, it is impossible for 
$x-y$ to be in $U^*\setminus Q$.
\end{remark}

\begin{corollary}
\label{co:std}
Suppose that $x$, $y \in \Phi$ are distinct. Then the
following are equivalent:
\begin{itemize}
\item [(i)] whenever $x \dom_W z \dom_W y$ for some $z \in \Phi$, then either
$z=x$ or $z=y$ (thus forming a cover of dominance);
\item [(ii)] there exists a $w\in W$ such that $wx \in D_0$ and $wy \in
-D_0$.  
\end{itemize}
\end{corollary}

\begin{proof}
Suppose that (i) is the case. 
Let $w$ be as in Proposition \ref{pp:key} above. First we show that then $wx\in
D_0$. Suppose for a contradiction that $wx \notin ~D_0$, and let $z \in D(wx)$. Then
Proposition~\ref{pp:key} yields that $wy\in \Phi^-$ and $(wy, z)\geq (wx, z) \geq 1 $. Hence
it is clear that $z \dom_W wy$. But this implies that $x \dom_W w^{-1}z \dom_W
y$ with $x \neq w^{-1}z \neq y$, contradicting (i). Therefore $wx \in D_0$, as
required. Exchanging the roles of $x$ and $-y$ we may deduce that $wy \in -D_0$.

Suppose that (ii) is the case and suppose for a contradiction that there exists
some $z\in
\Phi\setminus\{x, y\}$ such that $x\dom_W z\dom_W y$. Let $w\in W$ with $wx\in
D_0$ and $wy\in -D_0$. If $wz\in \Phi^+$ then
Lemma~\ref{lem:basicdom}~(ii) yields that $wx\dom_W wz$, contradicting the fact
that $wx\in D_0$. On the other hand, if $wz\in \Phi^-$, then
Lemma~\ref{lem:basicdom}~(ii) and (iii) yield that $-wy\dom_W -wz\in \Phi^+$,
contradicting the fact that $-wy\in D_0$.  

\end{proof}

Observe that applying Corollary \ref{co:std} to arbitrary reflection subgroup $W'$ of
$W$ yields the following:
\begin{corollary}
\label{co:W'}
Suppose that $W'$ is a reflection subgroup of $W$ with  $x$ and $y \in \Phi(W')$ being
distinct. Then the
following are equivalent:
\begin{itemize}
\item [(i)] whenever $x \dom_{W'} z \dom_{W'} y$ for some $z \in \Phi(W')$, then
either
$z=x$ or $z=y$;
\item [(ii)] there exists a $w\in W'$ such that $wx \in D_{W',\, 0}$ and $wy \in
-D_{W',\, 0}$.  
\end{itemize}
\qed
\end{corollary}

\begin{definition}
Suppose that $W'$ is a reflection subgroup of $W$ and $x, y\in \Phi(W')$ satisfy
both (i) and (ii) of Corollary
\ref{co:W'}. Then we say that the dominance between $x$ and $y$ is
\emph{minimal} with respect to $W'$.
\end{definition}

\begin{proposition}
\label{pp:mindom}
Suppose that $x,y \in \Phi$ are distinct with $x\dom_W y$, and let $W'$ be the 
dihedral reflection subgroup generated by $r_x$ and $r_y$. Then the dominance
between $x$ and $y$ with respect to $W'$ is minimal.

\end{proposition}

\begin{proof}
It follows from Corollary \ref{co:basic} (ii) that $x\dom_{W'} y$, and hence 
Lemma~\ref{lem:nodom} yields that $W'$ is infinite dihedral. Let 
$\Delta(W')=\{\,\alpha,\,\beta\,\}$. Then Proposition~\ref{pp:chain}~(i) yields that 
$D_{W',\, 0}=\{\,\alpha,\,\beta\,\}$.
 
On the other hand, it follows from Proposition~\ref{pp:candih} that there is some $w\in W'$ such
that either 
\begin{equation*}
 \left\{
\begin{array}{rl}
wx&=a\\
wy&=-b
\end{array} \right.\qquad\text{or else}\qquad
\left\{
\begin{array}{rl}
wx&=b\\
wy&=-a,
\end{array} \right.
\end{equation*}
consequently Corollary \ref{co:W'} yields that the dominance between $x$ and
$y$ with respect to $\langle\{ r_x, r_y \}\rangle$ is minimal.
\end{proof}
From the above proposition we may deduce:

\begin{proposition}
 Suppose that $x\in \Phi^+$ with $D(x)=\{x_1, x_2. \ldots, x_m\}$. For each
$i\in \{1, 2, \ldots, m\}$, set $W_i=\langle\{r_x, r_{x_i}\}\rangle$. Then
$W_i\neq W_j$ whenever $i\neq j$.
\end{proposition}

\begin{proof}
 For each $i\in \{1, 2, \ldots, m\}$, set $\{s_i, t_i\}=S(W_i)$. Suppose for a
contradiction that $W'=W_i=W_j$ for some $i\neq j$. Then we may write $\{s, t\}=\{s_i,
t_i\}=\{s_j, t_j\}$. Corollary~\ref{co:basic}~(ii) yields that $x\dom_{W_k}x_k$
for all $k\in \{1, 2, \ldots, m\}$, and since there is no non-trivial dominance
in finite Coxeter groups, it follows that $W_1, W_2, \ldots, W_m$ are all
infinite dihedral reflection subgroups. Hence it follows from Proposition~4.5.4
of \cite{ABFB} that $(\alpha_s, \alpha_t)\leq -1$. Set $c_n$ as in Proposition \ref{pp:rootsys}
for each $n\in \Z$.
Since $x\dom_W x_i$ and $x\dom_W x_j$, Proposition~\ref{pp:rootsys}~(ii) yields that either 
\begin{equation*}
 \left\{
\begin{array}{rl}
x &= c_{m}\alpha_s + c_{m+1}\alpha_t\\
x_i &= c_{m'}\alpha_s + c_{m'+1} \alpha_t\\
x_j &= c_{m''}\alpha_s+c_{m''+1}\alpha_t
\end{array} \right. \quad\text{or else}\quad
\left\{
\begin{array}{rl}
x &= c_{m}\alpha_s + c_{m-1}\alpha_t\\
x_i &= c_{m'}\alpha_s + c_{m'-1} \alpha_t\\
x_j &= c_{m''}\alpha_s+c_{m''-1}\alpha_t
\end{array} \right.
\end{equation*}
for some distinct integers $m, m'$ and $m''$. Observe that in either case 
$(x_i, x_j) \geq 1$,
 and therefore there will be (non-trivial)
dominance between $x_i$ and $x_j$. Without loss of any generality, we may assume
that $x\dom_W x_i\dom_W x_j$. Then 
$x\dom_{W'}x_i\dom_{W'} x_j$ by Corollary~\ref{co:basic}~(ii),
contradicting Proposition \ref{pp:mindom}.
\end{proof}

We close this paper with an alternative characterization for the imaginary cone $Q$ when
$\#\Delta<\infty$. 

\begin{proposition}
\label{pp:Q}
If $\#\Delta<\infty$ then
\begin{equation}
\label{eq:Qalt}
 Q=\{\, wv \mid w\in W\text{ and } v\in P \text{ such that } (v,a) \leq 0
\text{ for all } a\in \Phi^+ \,\}.
\end{equation}
\end{proposition}

\begin{proof}
First we denote the set on the right hand side of (\ref{eq:Qalt}) by $Z$, and
for each $b\in P$, define $\pos(b)=\{\, c\in \Phi^+\mid (b,c)>0\,\}$. Recall that
under the assumption that $\#\Delta<\infty$, Lemma~\ref{lem:tits} yields that
$$
Q=\{\,v\in \bigcap_{w\in W} wP\mid \text{$(v, a)\leq 0$ for all but finitely many $a\in\Phi^+$} \,\}.
$$

Let $u\in Q$ be arbitrary. Since $\#\Delta<\infty$, it follows from Lemma \ref{lem:tits} that $u\in P$.
If $\pos(u)=\emptyset$, then trivially $u\in Z$. 
Therefore we may assume that $\pos(u)\neq \emptyset$, and proceed by an induction on
$\#\pos(u)$ (this is only possible because $u\in Q$, and so $\#\pos(u)<\infty$).
Let $a\in \Delta$ be chosen such that $(u,a)>0$.  Then it can be readily checked that $\pos(r_a
u)=r_a(\pos(u)\setminus\{a\})$. Thus the inductive hypothesis yields that $r_a u
\in Z$. Clearly $Z$ is $W$-invariant, and so $u\in Z$, and hence $Q\subseteq Z$.

Conversely, if $x\in Z$, then $x=wv$ for some $w\in W$ and $v\in P$ such that
$(v,a)\leq 0$ for all $a\in \Delta$. Lemma~\ref{lem:l5} yields that $v\in U^*$, and since $U^*$
is clearly $W$-invariant, it follows that $x\in U^*$. Suppose that
$y\in \Phi^+$ with $(x, y)>0$. Since $(x,y) =(wv, y)=(v, w^{-1}y)$, and since
$(v, a)\leq 0$ for all $a\in \Phi^+$,  it follows
that $w^{-1}y \in \Phi^-$ and thus $y\in N(w^{-1})$. The finiteness of the set $N(w^{-1})$ then implies that 
$x\in Q$, and hence $Z\subseteq Q$.  
\end{proof}

\section{Acknowledgments}
A few results in Section \ref{sec:IMC} of this paper are taken from the author's PhD
thesis~\cite{FU0} and the author wishes to thank A/Prof.~R.~B.~Howlett for all
his help and encouragement throughout the author's PhD candidature. The author
also wishes to thank Prof.~G.~I.~Lehrer and Prof. ~R. Zhang for supporting this work.
Moreover, the author wishes to thank Prof. Dyer for his penetrating insight and valuable
suggestions. Finally, due gratitude must be paid to the referee of this paper, for her/his
extraordinary patience and an exhaustive list of helpful suggestions, in particular, we thank 
the referee for pointing to us the possibility of the converse of Theorem~\ref{th:cone}. 

\bibliographystyle{amsplain}

\begin{thebibliography}{4}
\bibitem{ABFB}
A. ~Bj\"{o}rner and F. ~Brenti, \emph{Combinatorics of Coxeter Groups}, Graduate
Texts in Mathematics, GTM 231, Springer, 2005 

\bibitem{BN68}
N.~Bourbaki, \emph{Groupes et algebras de Lie, Chapitres 4, 5 et 6 }, Hermann,
Paris, 1968


\bibitem{BH93}
B.~Brink and R.~B.~Howlett, \emph{A finiteness property and an automatic
structure of Coxeter groups}, Math. Ann. \textbf{296} (1993), 179--190.

\bibitem{BB94}
B.~Brink, \emph{On Root Systems and Automaticity of Coxeter Groups}, PhD thesis, 
University of Sydney. 1994.

\bibitem{BB98}
B.~Brink, \emph{The set of dominance-minimal roots}, J. Algebra \textbf{206}
(1998), 371--412.

\bibitem{VD82}
V. ~Deodhar, \emph{On the root system of a Coxeter group}, Comm. Algebra,
10(6):~611--630, 1982.

\bibitem{MD87}
M.~Dyer, \emph{Hecke algebras and reflections in Coxeter groups}, PhD thesis, 
University of Sydney, 1987.


\bibitem{MD90}
M. ~Dyer, \emph{Reflection Subgroups of Coxeter Systems}, J. Algebra
\textbf{135} (1990), 57--73.


\bibitem{bh}
M. ~Dyer, \emph{On the ``Bruhat graph'' of a Coxeter system},  Compositio Math. 
78  (1991),  no. 2, 185--191.

\bibitem{MD12}
M.~ Dyer, \emph{Imaginary cone and reflection subgroups of Coxeter groups}, Preprint
arXiv:1210.5206, http://arxiv.org/abs/1210.5206, October, 2012. 

\bibitem{TE09}
T. ~Edgar, \emph{Dominance and regularity in Coxeter groups}, PhD thesis,
University of Notre Dame, 2009
(downloadable from http://etd.nd.edu/ETD-db/).

\bibitem{FU0}
X. ~Fu, \emph{Root systems and reflection representations of Coxeter groups},
PhD thesis, University of Sydney, 2010.

\bibitem{FU1}
X. ~Fu, \emph{The dominance hierarchy in root systems of Coxeter groups},
J. Algebra, \textbf{366} (2012), 187--204.

\bibitem{HE1}
J.-Y. ~H\'ee, \emph{Le c\^one imaginaire d'une base de racines sur $\R$}, unpublished.

\bibitem{HE2}
J.-Y. ~H\'ee, \emph{Sur la torsion de Steinberg-Ree des groupes de 
Chevalley et des groupes de Kac-Moody}, PhD thesis, Universit\'e de Paris-Sud, Orsay, 1993. 


\bibitem{HH81}
H.~Hiller, \emph{Geometry of Coxeter Groups}. Research Notes in Mathematics 
54, Pitman (Advanced Publishing Program), Boston-London, 1981.

\bibitem{RB80}
R.~B.~Howlett, \emph{Normalizers of parabolic subgroups of reflection groups}, 
J. London Math. Soc. (2)  21  (1980), no. 1, 62--80.

\bibitem{RB96} 
R.~B.~Howlett, \emph{Introduction to Coxeter groups}, Lectures given at ANU,
1996 (available at http://www.maths.usyd.edu.au/res/Algebra/How/1997-6.html). 

\bibitem{HT97}
R.~B.~Howlett, P.~J.~Rowley and D.~E.~Taylor, 
\emph{On outer automorphisms of Coxeter groups},
Manuscripta Math. 94 (1997), 499--513. 

\bibitem{HM}
J. ~Humphreys, \emph{Reflection Groups and Coxeter Groups}, Cambridge studies in
advanced mathematics 29, 1990.

\bibitem{VK}
V.~G.~Kac, \emph{Infinite-dimensional Lie algebras}, third edition, Cambridge
University Press, Cambridge, 1990.


\bibitem{DK94}
D. ~Krammer, \emph{The conjugacy problem for Coxeter groups}, PhD thesis,
Universiteit Utrecht, 1994.

\bibitem{DK09}
D. ~Krammer, \emph{The conjugacy problem for Coxeter groups},  Groups Geom. Dyn., 
\textbf{3} (2009), No.1, 71--171.

\bibitem{SC91}
C.~Saunders, \emph{The Reflection Representations of Some Chevalley Groups}, J. Algebra
\textbf{137} (1991), 145--165.

\end{thebibliography}

\providecommand{\bysame}{\leavevmode\hbox to3em{\hrulefill}\thinspace}
\providecommand{\MR}{\relax\ifhmode\unskip\space\fi MR }
\providecommand{\MRhref}[2]{%
  \href{http://www.ams.org/mathscinet-getitem?mr=#1}{#2}
}
\providecommand{\href}[2]{#2}

\end{document}